\documentclass[11pt]{amsart}
\usepackage{anyfontsize}
\usepackage{tikz-cd}
\usepackage[headings]{fullpage}
\usepackage{amsmath}
\usepackage{amsfonts}
\usepackage{latexsym}
\usepackage{amssymb}
\usepackage{MnSymbol}
\newtheorem{thm}{Theorem}[section]

\newtheorem{lemma}{Lemma}[section]
\newtheorem{corr}{Corollary}[section]

\newtheorem{proposition}{Proposition}[section]

\def\ann{\mathop{\sf ann}}
\def\dep{\mathop{\sf depth}}
\def\dim{\mathop{\sf dim}}

\def\height{\mathop{\sf height}}

\def\sig{\sigma, \sigma^{-1}}
  
\begin{document}

\title[The Generic Degree of Autonomy]{The Generic Degree of Autonomy}

\maketitle

\begin{center}

\author{Shiva Shankar\footnote{Chennai Mathematical Institute,
Chennai (Madras), India} ~ and ~ Paula Rocha\footnote{Department of Electrical and Computer Engineering, University of Porto, Portugal}}

\vspace{.5cm}

{\em In memoriam: Professor Jan C. Willems (1939 - 2013)}

\end{center}

\begin{abstract} \noindent This paper calculates the degree of autonomy of a generic autonomous $n-D$ system defined by the kernel of a partial difference operator. The calculation implies that attaching a generic controller to a non-autonomous system, results in a controlled system whose degree of autonomy is the maximum possible. Thus a generic controller is maximally efficient with respect to the criterion of maximising degree of autonomy. The term  {\it generic} here refers to an open (dense) subset of the set of all $n-D$ systems with the Zariski topology. 
\end{abstract}

{\tiny \hspace{1cm}AMS classification: 39A14, 93B25, 13P25} \\

\section{Introduction} 

In this paper we study the notion of {\em degree of autonomy} of a multidimensional ($n-D$) autonomous system from Wood-Rogers-Owens \cite{wo} and Napp-Rocha \cite{d}. Autonomous systems arise when we attach a controller to a non-autonomous one, resulting in an autonomous system with prescribed properties. This procedure allows us to relate degree of autonomy to a notion of {\em strength} of a controller, and thence to a notion of its {\em efficiency}. Our discussion here is in the framework of J.C.Willems' Theory of Behaviors \cite{w}.
 
We study $n-D$ systems defined on the lattice $\mathbb{Z}^n$ whose laws are described by partial difference equations. A system is, by definition, a collection of functions $\{w: \mathbb{Z}^n \rightarrow \mathbb{C}^k\}$, these are all the possible `trajectories' of the system. The $k$ components of $w$ are the values of some $k$ attributes of the system; these values describe the system at various points of the lattice. A priori, any function could have been a trajectory, but the laws that govern the evolution of the system, proscribe certain trajectories from ever occurring. In the case of a linear shift-invariant system considered in this paper, these laws are linear constant coefficient partial difference equations, and we express the fact that a certain function is a possible trajectory, i.e. it satisfies these laws, by locating it in the kernel of a suitable operator on the space of all functions. Thus, a law limits the possible trajectories of the system, and other laws that the system satisfies cut down possibilities further. If the trajectories permitted by all the laws of the system need to be further restricted, then we construct a controller synthesised from additional laws, and attach it to the system forcing it to satisfy them as well, to arrive finally at  an acceptable collection of trajectories, specified  by an engineering situation. 

An autonomous system is one which, roughly speaking, does not admit {\em inputs} (formal definitions appear in Section 3). They are at the other end of the spectrum from controllable systems, which are non-autonomous systems with sufficiently many inputs that allows them to be steered from one trajectory to another. In Willems' worldview, the laws  a system obeys serve to reduce the availability of inputs to the system; this is the mechanism by which we reduce the system's possible trajectories as explained above. The imposition of sufficiently many laws, by the attachment of a controller, would eventually convert a non-autonomous system to an autonomous one, and it is this process that we study in this paper. 

We now explain another description of autonomous and non-autonomous linear systems that is important to the above interpretation of a controller. If a linear system satisfies some two laws, then it also satisfies every linear combination of these laws. Thus the set of all laws that a linear system satisfies has the structure of a module over the ring of difference equations. Necessarily then, a linear system satisfies an infinite number of laws. However as the ring of partial difference equations is Noetherian, it follows that the module of laws is finitely generated. If this module of laws can be generated by some $\ell$ number of laws, but not by any set of laws fewer than $\ell$ in number, then we say that the system `satisfies $\ell$ laws'.  Now, if this $\ell$ is less than the number $k$ of attributes of the system, then the system is said to be under-determined; otherwise it is over-determined.

The trajectories of a system are determined by its characteristic variety. This variety is all of $\mathbb{C}^n$ for a non-autonomous system, and a proper sub-variety of $\mathbb{C}^n$ for an autonomous one.  Thus, a non-autonomous system admits a very large number of possible trajectories, many of which might not meet various criteria, of boundedness, decay, stability and so on. A generic non-autonomous system is under-determined in the sense explained above, see \cite{s}, and its evolution is governed by fewer laws than the number of attributes, here $k$. By attaching a controller, that is by imposing additional laws the system must satisfy, we often construct an over-determined system, one which obeys more than $k$ laws, whose trajectories now satisfy the specified criteria (again generically, an over-determined system is autonomous). In this paper, we focus on one criterion, namely degree of autonomy.

Given an $n-D$ system, we can restrict it to an $m$ dimensional sub-lattice, i.e. to an embedding $\mathbb{Z}^m \hookrightarrow \mathbb{Z}^n$, to obtain an $m-D$ system. If the original system admitted inputs, so will the restricted system. However, it can be that the original system is autonomous, but its restriction to some embedding of $\mathbb{Z}^m$, for some $m$, is not. The degree of autonomy of a system measures this difference, it is defined to be the codimension of the largest sub-lattice such that the restriction of the system to it is not autonomous. Thus, the degree of autonomy of a non-autonomous system equals 0, and varies between 1 and $n$ for a nonzero autonomous system (a system is {\em strongly autonomous} when the degree equals $n$, \cite{lz,ps}). It is an important problem to calculate this invariant of a system, and the principal result in \cite{d,wo} states that it is equal to the codimension of the characteristic variety of the system.

To calculate the dimension of an affine variety is a difficult problem in general, but we show that it can be easily calculated for the characteristic varieties of an open subset of nonzero systems because we can calculate the length of a maximal regular sequence in their characteristic ideals.

Towards this, we first topologise the set of all $n-D$ systems, generalising the construction in \cite{s}. We then show that for a generic nonzero system, i.e. one belonging to a Zariski open (dense) set of systems, its characteristic ideal contains a regular sequence of this maximal length. As the ring of difference equations on $\mathbb{Z}^n$ is Cohen-Macaulay, we can then calculate the codimension of its characteristic variety, and hence its degree of autonomy. 

The above development leads to a notion of strength of a controller. Given a system, suppose the purpose of attaching a controller to it is to increase its degree of autonomy. This increase is the strength of the controller with respect to the system. The question then arises as to whether the increase in the degree of autonomy was accomplished by a controller synthesised from a minimum number of laws. This in turn leads to a notion of efficiency. It follows from the results of the paper that a generic controller is maximally efficient with respect to a given system, which is to say that a controller synthesised from a minimum number of {\em general} laws serves the stated purpose.

\noindent Remark: While the results in this paper are stated for systems defined over $\mathbb{Z}^n$, they are all equally valid for systems defined over $\mathbb{N}^n$ as well. This is because the principal result we use from \cite{wo} is valid in both cases.
  
\section{The topology of structured perturbations} 

We consider behavioral systems which arise as kernels of operators given by matrices with entries in the ring $A = \mathbb{C}[\sigma, \sigma^{-1}] := \mathbb{C}[\sigma_1, \sigma_1^{-1}, \ldots , \sigma_n, \sigma_n^{-1}]$ ($A$ is the Laurent polynomial ring, the coordinate ring of the complex $n$-torus $(\mathbb{C}^*)^n$). The term $\sigma_i$ acts on $w:\mathbb{Z}^n \rightarrow \mathbb{C}$ by shift of the $i$-th coordinate, namely $(\sigma_iw)(m_1, \ldots , m_n) = w(m_1, \ldots , m_i + 1, \ldots ,m_n)$. A monomial in $A$ acts on $w:\mathbb{Z}^n  \rightarrow \mathbb{C}$ by the  corresponding shift it defines and this action extends to $A$ by linearity: thus an element $a(\sigma, \sigma^{-1}) \in A$ defines an $A$-module map
$a(\sigma, \sigma^{-1}): (\mathbb{C})^{\mathbb{Z}^n}  \rightarrow  (\mathbb{C})^{\mathbb{Z}^n} $.
If $a(\sig) = (a_1(\sig), \ldots , a_k(\sig)) \in A^k$, then 
\[
\begin{array}{lccc}
a(\sig):  & (\mathbb{C}^k)^{\mathbb{Z}^n} &  \longrightarrow &  (\mathbb{C})^{\mathbb{Z}^n} \\
& w=(w_1, \ldots ,w_k) & \mapsto & \sum a_i(\sig)w_i
\end{array}
\]

By definition, the $n-D$ system defined by $a := a(\sig)$, or behavior $\mathcal{B}(a)$, is the subset of $(\mathbb{C}^k)^{\mathbb{Z}^n}$ given by the kernel of the above map. If $R$ is an $A$-submodule of $A^k$, then the behavior $\mathcal{B}(R)$ of $R$ is the intersection
\[ \bigcap_{a \in R} \mathcal{B}(a)
\]
If $R$ is generated by $\{a_1(\sig), \ldots ,a_\ell(\sig)\},~ a_i(\sig) = (a_{i1}(\sig), \ldots , a_{ik}(\sig))$, then $\mathcal{B}(R)$ is also the kernel of the operator 
\[
R(\sig): (\mathbb{C}^k)^{\mathbb{Z}^n}  \longrightarrow  (\mathbb{C}^\ell)^{\mathbb{Z}^n} 
\]
where $R(\sig)$ is the $\ell \times k$ matrix whose $(i,j)$-th entry is $a_{ij}(\sig)$. By the interpretation described in the introduction, the rows of 
$R(\sig)$ are the laws that govern the system, and to say that $w: \mathbb{Z}^n \rightarrow \mathbb{C}^k$ satisfies these laws is to say that it is in the kernel of the operator defined by it.

$\mathbb{C}^{\mathbb{Z}^n}$ is an injective $A$-module, and also a cogenerator \cite{o}. This implies that there is a bijective inclusion reversing correspondence between behaviors in  $(\mathbb{C}^k)^{\mathbb{Z}^n}$ and submodules of $A^k$.

Our principal purpose in this section is to topologize the set of all behaviors in $(\mathbb{C}^k)^{\mathbb{Z}^n}$, which by the above, is equivalent to topologizing the set of all submodules of $A^k$. We accomplish this by first topologizing the set of all matrices with $k$ columns and entries in $A$. While we broadly follow the procedure in \cite{s}, there are some differences. One difference is that now the ring $A$ is not the polynomial ring (i.e. the ring of constant coefficient partial differential operators) but its localization at $(\sigma_1 \sigma_2 ~\cdots ~\sigma_n)$. Another difference is that we do not fix the number of rows of the matrix, and this requires some comment.

We wish to topologize the set of matrices in order to study perturbations of a given behavior. An assumption made in \cite{s} is that the perturbations which are allowed are {\it structured} in such a way that the number of rows of a matrix does not change under perturbation.
We drop this assumption here, and allow changes in the number of rows as well. Thus we allow for the possibility that the number of laws that a behavior satisfies might itself change under perturbation. Such perturbations do arise in practice, for instance in the dynamics of switched systems, where sub-systems which were inactive are suddenly brought into play. Thus we generalise the notion of {\em structured perturbations} considered in \cite{s}. The question we study here, the maximum possible efficiency of controllers, can be calculated for a generic system, where the adjective `generic' is with respect to the more general notion of this paper.

In describing this topology of structured perturbations, we follow the procedure in \cite{s} with an additional step: matrices with $\ell_1$ rows are identified with matrices with $\ell_2$ rows, $\ell_2 > \ell_1$, by declaring all the entries in the last $(\ell_2-\ell_1)$ rows to equal 0. We topologize the set of all $\ell_2 \times k$ matrices in such a way that the set of $\ell_1 \times k$ matrices is a closed subspace, and then take the direct limit as the number of rows tends to infinity. 

In more detail,  let $k$ be a fixed positive integer, and let $\mathcal{M}_{\ell,k}$ be the set of matrices with $k$ columns and $\ell$ rows with entries in $A$. An element in $A$ is a sum of monomials $\sigma_1^{d_1}\cdots \sigma_n^{d_n}$ with complex coefficients, where the $d_i$ are integers (positive or negative). We define the {\em degree} of this monomial  to be $d = |d_1| + \cdots + |d_n|$, and the degree of an element in $A$ to be the maximum of the degrees of the monomial terms that it is the sum of.  Let $\mathcal{M}_{\ell,k}(d)$ be the set of those matrices in $\mathcal{M}_{\ell,k}$ whose entries are all bounded in degree by $d$. There are  

\begin{equation}
n_d = \sum_{j=0}^n 2^j {n \choose j}{d \choose j}
\end{equation}
 monomials
%\footnote{\phantom{xx}$^rC_s$ is the binomial coefficient $\frac{r!}{s!(r-s)!}$} 
in the $\sigma_i$ (i.e. with positive or negative powers) of degree bounded by $d$, there are $\ell k$ entires, hence an element $R(\sig)$ in  $\mathcal{M}_{\ell,k}(d)$ is a point in the affine space $\mathbb{C}^{N_d}$, equipped with the Zariski topology, where $N_d = n_d\cdot \ell k$ (and thus $N_d$ is a polynomial in $d$ of degree $n$). For $d_1 < d_2$, $\mathcal{M}_{\ell,k}(d_1)$ injects into $\mathcal{M}_{\ell,k}(d_2)$ as a closed subspace. Let the (strict) direct limit of the directed system of topological spaces $\{\mathcal{M}_{\ell,k}(d)\}_{d=0,1,2, \ldots }$ be also denoted $\mathcal{M}_{\ell,k}$.

As outlined above, we now consider the sequence $\{\mathcal{M}_{\ell,k}\}_{\ell = 1, 2, \ldots}$, where for $\ell_1 < \ell_2$, $\mathcal{M}_{\ell_1,k}$ injects into $\mathcal{M}_{\ell_2,k}$ as the closed subspace of matrices whose entries in the last $\ell_2 - \ell_1$ rows are all 0. This sequence is a directed system, and its (strict) direct limit is the space $\mathcal{M}(k)$ of all matrices with $k$ columns and entries in $A$. \\

\noindent Remark: Similar constructions involving the direct limit occur elsewhere in the literature, for instance \cite{nw}. Here, the authors topologise the direct limit using the euclidean topology on the finite dimensional spaces that define the direct system. In the algebraic setting of this paper, it is the Zariski topology that is relevant. As a consequence, our genericity statements are stronger than they would have been using the topology in \cite{nw}, for a Zariski open dense set is also open dense in the euclidean topology, and so is of full measure, of the second category, and so on. \\

Let $\mathcal{S}(k)$ be the set of all submodules of $A^k$. Consider the map
\[
\begin{array}{lccc}
\Pi_k: & \mathcal{M}(k) & \longrightarrow & \mathcal{S}(k) \\
  & R(\sig) & \mapsto & R
 \end{array}
\]
 where $\Pi_k$ maps a matrix with $k$ columns to the submodule of $A^k$ generated by its rows. We equip $\mathcal{S}(k)$ with the quotient topology (so that the map $\Pi_k$ is continuous). As the set $\mathcal{B}(k)$ of behaviors  in $(\mathbb{C}^k)^{\mathbb{Z}^n}$ is in bijective correspondence with elements in $\mathcal{S}(k)$, this procedure also makes $\mathcal{B}(k)$ a topological space. We refer to it as the Zariski topology on the set of behaviors in $(\mathbb{C}^k)^{\mathbb{Z}^n}$. A set of behaviors is said to be {\em generic} if it contains an open (dense) subset in this topology.

In particular, for $k=1$ we have a map
\[ \Pi_1: \mathcal{M}(1)  \longrightarrow  \mathcal{I}
\]
where we have denoted $\mathcal{S}(1)$, the space of ideals of $A$, by $\mathcal{I}$ ($\Pi_1$ maps a column in $\mathcal{M}(1)$ to the ideal generated by its entries). 

\begin{lemma} 
The map $\mathfrak{i}_k : \mathcal{S}(k) \longrightarrow \mathcal{I}$, mapping a submodule $R$ of $A^k$ to its characteristic ideal $\mathfrak{i}_k(R)$, i.e. to the 0-th Fitting ideal of $A^k/R$, is continuous (with respect to the Zariski topology described above).
\end{lemma}
\noindent Proof: Define a map $\frak{m}_k :\mathcal{M}(k) \rightarrow \mathcal{M}(1)$ by mapping $R(\sig)$ to its $k \times k$ minors (written in some fixed order). This map is continuous with respect to the Zariski topology described above as it involves only the arithmetic operations of addition and multiplication. 

Recall that the characteristic ideal of $R$, i.e. the 0-th Fitting ideal of $A^k/R$, is the ideal $\mathfrak{i}_k(R)$ generated by the $k \times k$ minors of any matrix whose rows generate $R$ (it follows from the Cauchy-Binet formula that this ideal is independent of the matrix whose rows generate $R$). Thus we have the commutative diagram
\begin{equation}
\begin{tikzcd}
\mathcal{M}(k) \arrow{r}{\Pi_k} \arrow[swap]{d}{\mathfrak{m}_k} & \mathcal{S}(k) \arrow{d}{\mathfrak{i}_k} \\
\mathcal{M}(1)  \arrow{r}{\Pi_1} & \mathcal{I}
\end{tikzcd}
\end{equation}
The maps $\mathfrak{m}_k$, $\Pi_k$ and $\Pi_1$ are continuous, and as $\mathcal{S}(k)$ has the quotient topology, $\mathfrak{i}_k$ is also continuous. 
\hspace*{\fill}$\square$\\

\noindent Remark: By means of the bijection between $\mathcal{B}(k)$ and $\mathcal{S}(k)$, the map $\frak{i}_k: \mathcal{B}(k) \rightarrow \mathcal{I}$, mapping a behavior in  $(\mathbb{C}^k)^{\mathbb{Z}^n}$ to its characteristic ideal, is continuous. \\

\noindent Remark on notation: When $\ell = k =1$, we denote $\mathcal{M}_{1,1}(d)$ by $A(d)$ and identify it with the affine space $\mathbb{C}^{n_d}$. Hence we denote the space $\mathcal{M}_{1,1}$ (the direct limit of the spaces $\{A(d)\}_{d=0,1, \ldots}$) also by $A$ and call its topology the Zariski topology. Consistent with this convention, the space $\mathcal{M}_{\ell,1}$ is denoted $A^\ell$.\\

We can restrict the above diagram to the closed subspace $\mathcal{M}_{\ell,k}$ of $\mathcal{M}(k)$: let $\mathcal{S}_{\ell,k}$ be the set of submodules of $A^k$ that can be generated by $\ell$ elements, and let $\mathcal{B}_{\ell,k}$ be the set of $n-D$ behaviors in bijective correspondence with it. We denote $\mathcal{S}_{\ell,1}$ by $\mathcal{I}_\ell$ (thus $\mathcal{I}_\ell$ is the set of ideals of $A$ that can be generated by $\ell$ elements). Then the map $\Pi_k$ restricts to a surjection $\Pi_{\ell,k} :\mathcal{M}_{\ell,k} \rightarrow \mathcal{S}_{\ell,k}$. The quotient topology on $\mathcal{S}_{\ell,k}$ makes it a closed subspace of $\mathcal{S}(k)$, and hence $\mathcal{B}_{\ell,k}$ a closed subspace of $\mathcal{B}(k)$. The map $\frak{m}_k$ restricts to $\frak{m}_{\ell,k}: \mathcal{M}_{\ell,k} \rightarrow 
A^{\small{\ell \choose k}}$, 
mapping an $\ell \times k$ matrix to its $\ell \choose k$
many $k\times k$ minors. Similarly $\frak{i}_k$ restricts to a map $\frak{i}_{\ell,k}$ to give the commutative diagram
\begin{equation}
\begin{tikzcd}
\mathcal{M}_{\ell,k} \arrow{r}{\Pi_{\ell,k}} \arrow[swap]{d}{\mathfrak{m}_{\ell,k}} & \mathcal{S}_{\ell,k} \arrow{d}{\mathfrak{i}_{\ell,k}} \\
\phantom{xx} A^{\tiny{\ell \choose k}}
\phantom{xx} \arrow{r}{\Pi_{\tiny{\ell \choose k},1}} & ~\phantom{x} \mathcal{I}_{\tiny{\ell \choose k}}
\end{tikzcd}
\end{equation}
This implies as before that the map $\frak{i}_{\ell,k}: \mathcal{S}_{\ell,k} \rightarrow \mathcal{I}_{\small{\ell \choose k}}$, and the corresponding map $\mathcal{B}_{\ell,k} \rightarrow \mathcal{I}_{\small{\ell \choose k}}$,
are continuous. (We follow the usual convention that ${\ell \choose k} = 0$ 
for $\ell< k$, and that the ideal generated by the empty set is the 0 ideal.)\\

\noindent Notation: We denote by $A_+$ the subring $\mathbb{C}[\sigma] := \mathbb{C}[\sigma_1, \ldots ,\sigma_n]$ of $A$, and by $A_+(d)$ the set of polynomials of degree bounded by $d$. There are ${n+d \choose n}$ monomials in the $\sigma_i$ with non-negative exponents and degree bounded by $d$, hence  $A_+(d)$ can be identified with the affine space $\mathbb{C}^{\tiny {n+d \choose n}}$ (with the Zariski topology). The direct limit of the spaces $\{A_+(d)\}_{d=0,1, \ldots}$ is the space $A_+$.  

The ring $A$ is the localisation of $A_+$ at the multiplicative set $\{(\sigma_1 \sigma_2 \cdots \sigma_n)^j ~| ~j \geq 0\}$ generated by the product of the $\sigma_i$. As $(\sigma_1 \sigma_2 ~\cdots ~\sigma_n)^{-1}$ exists in $A$, so also do the $\sigma_i^{-1}$, $i = 1, \ldots n$. Thus the units in $A$ are the nonzero constants (the units of $A_+$) together with all the monomials $\sigma_1^{d_1} \sigma_2^{d_2} \cdots \sigma_n^{d_n}$, $d_i \in \mathbb{Z}$, in the $\sigma_i$ (with nonzero complex coefficients).

An element of the Laurent polynomial ring $A$ is a sum of monomial terms $\sigma ^d := \sigma_1^{d_1} \cdots \sigma_n^{d_n}$ with complex coefficients; in the topological space $A$ it corresponds to the point whose coordinates are these coefficients (if a monomial does not appear in the sum, then its coefficient is 0). Thus, there is a `coordinate axis' in the space $A$ corresponding to each monomial in the ring $A$. We denote the indeterminate, and the coordinate axis, corresponding to the monomial $\sigma ^d$ by $X_d$ or $X_{ d_1 \cdots d_n}$. The points of the space $A$ corresponding to the units in the ring $A$ are the points on the axes minus the origin. Its closure is the union of the coordinate axes, and is a proper Zariski closed subset of the space $A$.

$A_+$ is the subspace of $A$ given by points whose $X_{d_1 \cdots d_n}$-coordinates equal 0 whenever any of the exponents $d_i$ is negative. It is a Zariski closed subspace of $A$.\\

The {\em spaces} $A$ and $A_+$ are not Noetherian, for instance the descending sequence of Zariski closed subspaces $A_+ \supsetneq \{X_1 = 0\} \supsetneq \{X_1 = 0, X_2 = 0\} \supsetneq ...$ ~does not stabilise. If $X$ denotes the union of the indeterminates $X_d$, and $X_+$ denotes the union of the $X_d$ with $d = (d_1, \ldots ,d_n)$ such that the $d_i$ are all non-negative, then the coordinate rings of the spaces $A$, $A_+$ are $\mathbb{C}[X]$ and $\mathbb{C}[X_+]$ respectively. 

We collect a few elementary properties of the spaces $A$ and $A_+$.

\begin{lemma}
\noindent (i)  The spaces $A$ and $A_+$ are irreducible. \\
\noindent (ii) The projections $A^{r+s} \rightarrow A^r$ and $A_+^{r+s} \rightarrow A_+^r$ are open.
\end{lemma}

\noindent Proof: These statements follow from corresponding statements for $A(d)$ and $A_+(d)$.
 \hspace*{\fill}$\square$\\
 
By (i) every nonempty open subset of $A$ or $A_+$ is dense. Nonetheless, we sometimes use
the phrase `open dense' for emphasis. \\
 
\noindent Remark: As we have already remarked, all the results of this paper also hold for linear systems in $(\mathbb{C}^k)^{\mathbb{N}^n}$ defined by subodules of $A_+^k$ although we do not explicitely mention this fact hereafter.

\section{Degree of autonomy of a generic system} 

We recollect from the introduction the notion of degree of autonomy of an $n-D$ system, introduced in \cite{wo} and studied further in \cite{d}. \\

\noindent Definition: Let $\mathcal{B}(R) \subset (\mathbb{C}^k)^{\mathbb{Z}^n}$ be the $n-D$ system defined by the submodule $R$ of $A^k$. Then $\mathcal{B}(R)$ is autonomous if none of the projections $\pi_i: \mathcal{B}(R) \rightarrow \mathbb{C}^{\mathbb{Z}^n}$, $(w_1, \ldots , w_k) \mapsto w_i$, $i = 1, \ldots ,k$, is surjective.

The characterisation in \cite{s} of autonomous systems defined in the space of distributions or smooth functions by partial differential equations, carries over to the case of $n-D$ systems.\\

\noindent Proposition \cite{s}: The $n-D$ system $\mathcal{B}(R)$ is autonomous if and only if the annihilator $\ann(A^k/R)$ is nonzero.
\hspace*{\fill}$\square$\\

\noindent Remark: The characteristic ideal of the submodule $R$ and the annihilator of $A^k/R$ are related, for instance \cite{s}:

\hspace{5cm} $\frak{i}_k(R) \subseteq \ann(A^k/R) \subseteq \sqrt{\frak{i}_k(R)}$ 

\noindent Thus the condition $\ann(A^k/R) \neq 0$ is equivalent to the condition that $\frak{i}_k(R) \neq 0$. An autonomous system is therefore also defined by the condition that its characteristic ideal is nonzero, for instance \cite{ps}. It follows that an autonomous system is necessarily over-determined (conversely, a generic over-determined system is autonomous \cite{s}). 

These conditions are equivalent to saying that $A^k/R$ is a torsion module. This implies that the behavior of $R$ does not contain any nonzero controllable sub-behaviors, \cite{s}. \\

A choice of $m$ indices $1 \leq i_1 < i_2 < \ldots < i_m \leq n$ defines an $m$ dimensional sub-lattice of $\mathbb{Z}^n$ as follows: it is the inclusion $\iota: \mathbb{Z}^m \hookrightarrow \mathbb{Z}^n$ mapping $z=(z_1, \ldots , z_m)$ to $z'=(z'_1, \ldots , z'_n)$ where $z'_{i_j} = z_j, 1 \leq j \leq m$, and the other coordinates equal to 0. These indices also define the subring $A_\iota =  \mathbb{C}[\sigma_{i_1}, \sigma_{i_1} ^{-1}, \ldots , \sigma_{i_m}, \sigma_{i_m}^{-1}]$ of $A$.

Let $B = \{w:\mathbb{Z}^n \rightarrow \mathbb{C}^k \}$ be an $n-D$ system, and let $\iota : \mathbb{Z}^m \hookrightarrow \mathbb{Z}^n$ be a sub-lattice. The restriction of $B$ to $\mathbb{Z}^m$ is the $m-D$ system $\iota^* (B) = \{w \circ \iota: \mathbb{Z}^m \rightarrow \mathbb{C}^k ~| ~w \in B \}$. If $B$ is defined by the submodule $R$ of $A^k$, then $\iota^* (B)$ is the behavior defined by the $A_\iota$-submodule $\iota^* R = R \cap A_\iota^k$, \cite{d,wo}. As $\ann(A_\iota^k/\iota^*R) = \ann(A^k/R) \cap A_\iota$, it follows (from the proposition quoted above) that if $B$ is not autonomous, so also is $\iota^* (B)$ not autonomous. However it may be that $B$ is autonomous but $\iota^* (B)$ is not, for some sub-lattice $\iota$. This motivates the following definition. \\

\noindent Definition \cite{d,wo}: Let $B$ be an $n-D$ system.
Let $m \geq 0$ be the largest integer such that for some sub-lattice $\iota : \mathbb{Z}^m \hookrightarrow \mathbb{Z}^n$, the restricted 
$m-D$ system $\iota^* (B)$ is not autonomous. Then the degree of autonomy of $B$ equals $n-m$. If $B = 0$, then its degree of autonomy is 
defined to be $\infty$. \\ 

Thus the degree of autonomy of a non-autonomous system equals 0. \\

The affine variety in $\mathbb{C}^n$ defined by the characteristic ideal of a system is its characteristic variety.
A principal result in Wood et al. \cite{wo} and Napp-Rocha \cite{d} is the relationship between the degree of autonomy of a system and the dimension of its characteristic variety. \\

\noindent  {\bf Theorem} \cite{d,wo}: The degree of autonomy of an $n-D$ system equals the codimension of its characteristic variety. 
\hspace*{\fill}$\square$ \\ 

The characteristic variety of a strongly autonomous $n-D$ system is a finite set of points in $\mathbb{C}^n$ (by definition), hence its degree of autonomy  equals $n$, \cite{lz,ps}. \\

\noindent Remark: Suppose $B_1 = \mathcal{B}(R_1)$ is a sub-system of $B = \mathcal{B}(R)$, then $R \subset R_1$. By the above theorem, the degree of autonomy of $B_1$ is greater than or equal to that of $B$. Suppose this increase in the degree of autonomy was the result of attaching a controller $C = \mathcal{B}(R')$ to $B$, then $R_1 = R + R'$. The purpose of this paper is to study the following question: what is the smallest $C$ (i.e. one defined by as small a submodule $R'$ of $A^k$ as possible) that results in a given increase in the degree of autonomy of $B$? \\

While generators of the characteristic ideal of a system defined by a submodule $R$ of $A^k$ can be easily calculated $-$ they are the $k \times k$ minors of any matrix whose rows generate $R$ $-$ the calculation of the dimension of its variety, and hence the degree of autonomy of $\mathcal{B}(R)$, is a very difficult problem in general (there are now computer packages based on Gr\"obner basis methods that calculate the dimension of a variety). The rest of this section is about this problem, and is organised as follows: we first isolate a large class of ideals where this calculation is indeed possible. More precisely, we show in Proposition 3.4 below that points in $A^r$, $r \leq n$, corresponding to regular sequences, are open dense in $A^r$. The varieties of these ideals are of codimension $n-r$. We then show that the set of nonzero  $n-D$ systems whose characteristic ideals contain regular sequences of maximum length, is open in the set of all nonzero systems. This allows us to calculate the degree of autonomy for this class of open systems. 

We illustrate the nature of the above problem with an elementary example. \\

\noindent {\bf Example} ($2-D$ scalar systems):  Now $A = \mathbb{C}[\sigma_1, \sigma_1^{-1}, \sigma_2, \sigma_2^{-1}]$ and $k = 1$; scalar systems are thus defined by ideals of $A$. While there is no bound on the number of elements needed to generate ideals of $A$, generically only the cases described below occur: \\

\noindent (i) $\ell = 1$, i.e. the system is defined by a principal ideal $(a(\sig))$. The system is then given by the kernel of

\hspace{5cm} $a(\sig): (\mathbb{C})^{\mathbb{Z}^n} \rightarrow (\mathbb{C})^{\mathbb{Z}^n}$

By Krull's Principal Ideal Theorem \cite{k}, a (proper) principal ideal defines a codimension 1 variety in $\mathbb{C}^2$, hence the degree of autonomy of this system equals 1. The set of elements of $A$ that are not units is Zariski open in the space $A$, hence there is an open dense set of $2-D$ systems in $\mathcal{B}_{1,1}$ whose degree of autonomy equals 1. \\

\noindent (ii) $\ell = 2$, i.e. the system is defined by an ideal that can be generated by 2 elements (the previous case of $\ell = 1$ is therefore included here). The system is given by the kernel of a map 
\[ \left(\begin{array}{c}a_1(\sig) \\a_2(\sig) \end{array} \right): (\mathbb{C})^{\mathbb{Z}^n} \rightarrow (\mathbb{C}^2)^{\mathbb{Z}^n}
\] By Proposition 3.6 below, the set of elements $(a_1, a_2) \in A^2$ such that $a_1 \neq a_2$ and both irreducible, is open dense in $A^2$. The variety of the ideal generated by such a pair is a finite set of points in $\mathbb{C}^2$, and the degree of autonomy of such systems, an open dense subset of $\mathcal{B}_{2,1}$, equals 2. As $\mathcal{B}_{1,1}$ is Zariski closed in $\mathcal{B}_{2,1}$, systems whose degree of autonomy equals 2 are also open dense in $\mathcal{B}_{2,1} \setminus \mathcal{B}_{1,1}$. \\

\noindent (iii) $\ell \geq 3$, i.e. the system is defined by an ideal that can be generated by 3 elements. By Corollary 3.1 below, an open dense set of points in $A^{\ell}$ define the unit ideal. The system in $\mathcal{B}_{\ell,1} \setminus \mathcal{B}_{2,1}$ corresponding to these points is the zero behavior, whose degree of autonomy equals $\infty$.
\hspace*{\fill}$\square$\\

In the rest of this section we extend this example to general $n-D$ systems (Theorem 3.2 below).

\begin{lemma}
Let $\mu: A^r \times A^r \rightarrow A$ be the map $((a_i), (b_i)) \mapsto \sum_1^r a_ib_i$. Then $\mu ^{-1} (0)$ is a proper Zariski closed subset of the space $A^r \times A^r$.
\end{lemma}
 
 \noindent Proof: The map $\mu$ restricts to $\mu _d: A(d)^r \times A(d)^r \rightarrow A(2d)$. It is given by adding and multiplying the coefficients of the Laurent polynomials $a_i$ and $b_i$, hence it is continuous in the Zariski topology. The point 0 is closed in $A$, hence $\mu_d ^{-1} (0)$ is closed in $A(d)^r \times A(d)^r$. The direct limit of $\{\mu_d^{-1}(0)\}_{d=0,1,\ldots}$ equals $\mu^{-1}(0)$, hence it is closed in $A^r \times A^r$.
 \hspace*{\fill}$\square$ \\
 
 We need to study the image of the map $\mu$.
 
\begin{lemma}
Let $I$ be an ideal of $A$. Then $I$ is a Zariski closed subset of the space $A$. It is a proper closed subset if and only if $I$ is a proper ideal; then the set of elements of $A$ which are not in $I$ is an open dense subset of the space $A$. 
\end{lemma}
\noindent Proof: Let $I$ be generated by $\{a_1, \ldots ,a_r\}$, and let $a$ denote the point $(a_1 \ldots a_r) \in A^r$ . Define the map $\mu_a: A^r \rightarrow A$ by $\mu_a(b_1, \ldots ,b_r) = \sum a_ib_i$. This is a $\mathbb{C}$-linear map as the coefficients of the Laurent polynomial $\sum a_ib_i$ are $\mathbb{C}$-linear combinations of the coefficients of the $b_i$. Its image is precisely the ideal $I$, and hence $\mu_a$ is not surjective exactly when $I$ is a proper ideal of $A$. 

Assume then that $I$ is proper. Let the maximum of the degrees of the $a_i$ be $s$. For each $d$, the map $\mu_a$ restricts to a map $\mu_a(d): (A(d))^r \rightarrow A(d+s)$. Its image, say $I(d + s)$, is contained in $I \cap A(d + s)$.
As $I$ is proper, $I(d + s)$ is a proper linear subspace of $A(d+s) (\simeq \mathbb{C}^{n_{d+s}}$). It is therefore a proper Zariski closed subset of $A(d+s)$ whose vanishing ideal is generated by linear forms (in the indeterminates $\{X_d\}$). Its complement is then open dense. 
For $d_1 < d_2$, the map $\mu_a(d_2)$ restricts to $\mu_a(d_1)$, hence $I(d_1 + s)  \subset I(d_2 + s)$. The direct limit of these closed subspaces $\{I(d + s)\}_{d = 0, 1, \ldots}$ is $I$, hence $I$ is a proper Zariski closed subset of $A$.  
\hspace*{\fill}$\square$\\

\noindent Remark:  More generally, every $\mathbb{C}$-linear subspace of $A$ is Zariski closed, given by the common zeros of linear forms in the indeterminates $\{X_d\}$. So is therefore every affine linear subset of $A$. \\

\noindent Remark: If $I = (a)$ is a nonzero principal ideal, generated by a Laurent polynomial of degree $s$, then the map $\mu_a(d): A(d) \rightarrow A(d+s)$ is injective ($A$ is an integral domain).  Its image is therefore a $\mathbb{C}$-subspace of codimension $n_{d+s} - n_d$ (notation as in equation (1)), an increasing function of $s$.  If $a$ is in $A_+$ (i.e. $a$ is a polynomial), then the map $\mu_a(d)$ restricts to $\mu_a(d): A_+(d) \rightarrow A_+(d+s)$. Its image is of codimension ${n+d+s \choose n} - {n+d \choose n}$, also an increasing function of $s$.

\begin{proposition}
For $a = (a_1, \ldots , a_r) \in A^r$, let $I_a$ be the Zariski closed subset of $A$ consisting of the elements in the ideal $(a_1, \ldots ,a_r)$ generated by the coordinates $a_1, \ldots , a_r$ of $a$ $($namely, the above lemma$)$. Then the set $\frak{I} = \{(a, a_{r+1}) ~| ~ a \in A^r, ~ a_{r+1} \in I_a\}$ is Zariski closed in the space $A^r \times A$.
\end{proposition}

\noindent Proof: It suffices to observe that in the proof of Lemma 3.2 above, $\mu_a$ is a polynomial function of $a = (a_1, \ldots, a_r)$, i.e. in the coefficients of $a_1, \ldots , a_r$ as well, namely Lemma 3.1. Indeed, $\mu_a$ is linear in the coefficients of the $a_i$, and thus the vanishing ideal of  $\frak{I}$ is generated by 2-forms (in the indeterminates $\{X_d\}$).
\hspace*{\fill}$\square$\\

\begin{lemma}
Let $I$ be a proper ideal of $A$. Then the set of zero divisors on $A/I$ is closed in the space $A$, and hence the set of nonzero divisors on 
$A/I$ is open dense.
\end{lemma}

\noindent Proof: The set of zero divisors on $A/I$ is the union of its finite number of associated primes, and this finite union is closed in $A$ (by Lemma 3.2).
\hspace*{\fill}$\square$\\

Let $B = \mathbb{C}[\sigma_1 - \sigma_1^{-1}, \ldots , \sigma_n - \sigma_n^{-1}]$. The $(\sigma_i - \sigma_i^{-1}), 1 \leq i \leq n,$ are algebraically independent, and $A$ is integral over B. 
As $B$ has Krull dimension $n$, so does $A$. \\

Suppose $I$ is a maximal ideal in $A$, it is then generated by $n$ elements, say $I = (a_1, \ldots ,a_n)$.  The sum $I + (a)$ is not equal to $A$ if and only if the element $a$ belongs to $I$. Thus $I + (a) = A$ for $a$ in a dense open subset of $A$. Clearly this is equivalent to the statement that $A/I \simeq\mathbb{C}$, for in $\mathbb{C}$ the set of units is Zariski open.

Further, let $M$ be the set of points $a = (a_1, \ldots ,a_n) \in A^n$ such that the ideal $(a_1, \ldots ,a_n)$ generated by its coordinates is maximal (we show below that $M$ is open dense in the space $A^n$). Then by Proposition 3.1, the set of points $(a_1, \ldots , a_n, a_{n+1}) \in A^{n+1}$ such that $(a_1, \ldots ,a_n) \in M$ and the ideal generated by these $n+1$ coordinates is proper, is closed in the space $M \times A$. 

The opposite is however the case for proper ideals generated by fewer than $n$ elements, namely Proposition 3.2 below, and is suggested by the following heuristic:

Let $I$ be a proper ideal generated by $r < n$ elements, then the height of $I$ is at most $r$ (by Krull's Height Theorem), hence the dimension of $A/I$ is at least $n-r$. By Noether normalization (for instance \cite{k}), $A/I$ is isomorphic to an integral extension of a polynomial ring with number of indeterminates at least $n-r$. Therefore, the set of units in $A/I$ is contained in a proper Zariski closed set, and the set of elements $a \in A$ such that the sum $I + (a)$ is not equal to $A$ contains an open dense subset of the space $A$.

\begin{proposition}
Let $r \leq n$. Then the set $P_r$ of elements $(a_1, \ldots ,a_r) \in A^r$ such that the ideal $(a_1, \ldots ,a_r)$ generated by its coordinates is a proper ideal of $A$, contains an open dense subset of the space $A^r$.
\end{proposition}

We use the following result of Brownawell \cite{b}:

\noindent Theorem (Brownawell): Suppose the ideal generated by the polynomials $p_1, \ldots ,p_m$ equals $A_+ = \mathbb{C}[\sigma_1, \ldots, \sigma_n]$, where the degree of the $p_i$ is less than or equal to $d$. Then there are polynomials $q_1, \ldots ,q_m$ such that $\sum p_i q_i = 1$, where the degree of the $q_i$ is less than or equal to $D = n^2 d^n + nd$. \\

\noindent Proof of proposition:  It suffices to prove the statement for $r = n$, for by Lemma 2.2 (ii), if $P_n$ contains an open dense subset of $A^n$, then its projection to $A^r$ also contains an open dense subset. The coordinates of these points generate proper ideals of $A$. Thus this projection of $P_n$ is contained in $P_r$.

We first prove the statement for $A_+$, i.e. we prove that the set of points in $A_+^n$ whose coordinates generate proper ideals of $A_+$, contains an open dense subset of $A_+^n$.

So let $I$ be the ideal generated by the coordinates of a point $a = (a_1, \ldots ,a_n) \in A_+^n$. 
 Let the degrees of the $a_i$ be bounded by $d$. To say that $I \subsetneq A_+$ is to say that  
the map $\frak{a}: \mathbb{C}^n \rightarrow \mathbb{C}^n$ defined by $x \mapsto (a_1(x), \ldots, a_n(x))$ has nonempty inverse image $\frak{a}^{-1}(0)$. Clearly there is a point $a$ in $A_+(d)^n$ such that the corresponding map $\frak{a}$ has a regular point $x$ in this inverse image. This means that the rank of $\frak{a}$ equals $n$ at $x$, hence all smooth maps sufficiently close to $\frak{a}$ in the compact-open topology also include 0 in their images (Inverse Function Theorem). Restricting to maps given by elements in $A_+(d)^n$ as above, this means that there is an open neighbourhood of  $a$ in the {\em euclidean topology} on $A_+(d)^n$ such that the variety of the ideal generated by the coordinates of every point in it, is nonempty. These ideals are thus proper ideals of $A_+$.

Thus $P_n \cap A_+(d)^n$ contains a euclidean open subset of $A_+(d)^n$, and we now show that it contains a nonempty Zariski open subset.

Suppose to the contrary that $I$ (in the above notation) were not a proper ideal of $A_+$. 
Then by the theorem of Brownawell, there are polynomials $b_1, \ldots ,b_n$ of degree bounded by $D = n^2 d^n + nd$, such that $\sum a_i b_i = 1$. This is a Zariski closed condition on the coefficients of the $a_i$ and $b_i$ (as in Lemma 3.1), and it defines a proper affine variety in the affine space $A_+(d)^n \times A_+(D)^n$. Its projection $C$ to the first $n$ coordinates in $A_+(d)^n$ is a constructible set (by Chevalley's theorem (EGA IV, 1.8.4) quoted in \cite{k}, the image of a variety is a constructible set). The coordinates of a point in $C$ generates the unit ideal $A_+$ (hence the unit ideal of $A$), and the coordinates of points in the complement $C' = A_+(d)^n \setminus C$, which is also constructible, generate proper ideals of $A_+$. 

We have shown at the outset that $C'$ contains a euclidean open subset; as it is constructible, it must therefore contain a nonempty Zariski open subset, say $O_+(d)$, of $A_+(d)^n$. This is true for every $d$, hence the set of points in $A_+^n$ which generate a proper ideal of $A_+$, contains an open dense subset of $A_+^n$.

Every ideal of $A$ is an extension of an ideal of $A_+$, and the proper ideals of $A_+$ that do not remain proper in $A$ are those whose varieties (in $\mathbb{C}^n$) are contained in the varieties of the $\sigma_i$ (that is, those varieties contained in the codimension 1 coordinate hyper-planes $\sigma_i = 0$). This is a Zariski closed condition, for let $d > 0$, and let $s=\sigma_1^d
 \sigma_2^d \cdots \sigma_n^d$. Then $s$ defines by multiplication an injective map, $s:A(d)^n \rightarrow A(d+nd)^n$, where $(a_1, \ldots a_n) \mapsto (sa_1, \ldots ,sa_n)$. Its image is contained in $A_+(d+nd)^n$, and the inverse image $O(d)$, of  the open set $O_+(d+nd)$, is open  in $A(d)^n$. The set of points corresponding to proper ideals in $A_+$ but which extend to the full ring $A$ is Zariski closed in $O(d)$, and its complement is contained in $P_n \cap A(d)^n$. 
 \hspace*{\fill}$\square$\\

\begin{proposition}
Let $r < n$. The set $U_{r+1}$  of points $(a_1, ..., a_{r+1}) \in P_{r+1}$ such that $a_1 \neq 0$ in $A$, $a_2 \neq 0$ in $A/(a_1)$, ..., $a_{r+1} \neq 0$ in $A/(a_1, \ldots ,a_r)$, is open dense in $P_{r+1}$, and hence contains an open dense subset of $A^{r+1}$.
\end{proposition}

\noindent Proof:  The point 0 is closed in $A$, hence the statement is true for $r = 0$. Assume by induction that $U_r$ is open dense in $P_r$. The subset $\{(a,b_{r+1}) ~| ~a \in U_r, ~b_{r+1} \in I_a\}$ of the closed set $\frak{I}$ of Proposition 3.1, is closed in $U_r \times A$. Its complement is open dense in $U_r \times A$; hence $U_{r+1}$, which is the intersection of this complement with $P_{r+1}$, is open dense in $P_{r+1}$, and hence contains a dense open subset of $A^{r+1}$.   
\hspace*{\fill}$\square$\\

\begin{proposition}
Let $r < n$. The set $N_{r+1}$ of points $(a_1, \ldots , a_{r+1}) \in U_{r+1}$ such that $a_1$ is a nonzero divisor $($nzd$\phantom{`})$ on $A$, $a_2$ is a nzd on $A/(a_1)$, ..., $a_{r+1}$ is a nzd on $A/(a_1, \ldots ,a_r)$ is open dense in $U_{r+1}$, and hence contains an open dense subset of $A^{r+1}$. 
\end{proposition}

\noindent Proof: The statement is true for $r = 0$ as now $N_1$ equals $U_1$ ($A$ is an integral domain). Assume by induction that  $N_r$ is open dense in $U_r$. For $a = (a_1, \ldots ,a_r)$ in $N_r$, let $Z_a$ be the set of zero divisors on $A/(a_1, \ldots ,a_r)$. $Z_a$ is a closed subset of the space $A$ by Lemma 3.3. We need to show that the set $\frak{Z} = \{(a, b_{r+1}\} ~| ~a \in N_r, b_{r+1} \in Z_a\}$ is closed in $N_r \times A$, for its complement in $U_{r+1}$ is precisely $N_{r+1}$. 
To show this, it suffices by Proposition 3.1 to show that every element in the vanishing ideal of $Z_a$ is a polynomial function of elements in the vanishing ideal of the closed set $I_a$ of the space $A$.

Let $\mu: A \times A \rightarrow A$ be the multiplication map, mapping $(a, b)$ to $a b$ (the map of Lemma 3.1 for $r = 1$). It is an algebraic map as it is given by multiplying and adding various coefficients. Denote also by $\mu$ its restriction $\mu : A \times (A \setminus I_a) \rightarrow A$. Then $V = \mu ^{-1}(I_a)$ is a Zariski closed subset of the space $A \times (A \setminus I_a)$. Let its vanishing ideal be $J$; its elements are polynomial functions of the elements of the vanishing ideal of $I_a$, and therefore polynomial functions of the coefficients of the components $a_1, \ldots , a_r$ of $a$. 

The set $Z_a$ is the projection of $V$ to the first factor $A$. In general a projection is not Zariski closed as the space $A$ is not complete, but here $Z_a$ is indeed closed (Lemma 3.3). Hence its vanishing ideal equals $\iota ^{-1} (J)$, where $\iota$ is the inclusion of $A$ in $A \times (A \setminus I_a)$.  

As $\iota$ is an algebraic map, it follows that elements of  $\iota ^{-1} (J)$ are polynomial functions of the coefficients of $a_1, \ldots ,a_r$. This completes the proof. 
\hspace*{\fill}$\square$\\

Recall the definition of a Cohen-Macaulay ring:  a sequence $a_1, \ldots ,a_t$ in a Noetherian ring $R$ is regular if (i) the ideal $(a_1, \ldots , a_t) \subsetneq R$, and (ii) $a_1$ is a nonzero divisor in $R$, and for each $i$, $2 \le i \le t$, $a_i$ is a nonzero divisor on $R/(a_1, \dots , a_{i-1})$ (thus Proposition 3.4 asserts that for $r < n$, the set of points in $A^{r+1}$ corresponding to regular sequences contains an open dense subset of $A^{r+1}$). The depth of an ideal $I$ is the length of any maximal regular sequence in $I$. The ring $R$ is Cohen-Macaulay if for every ideal $I$ of $R$, $\dep(I) = \height(I)$. The localisation of a Cohen Macaulay ring at any multiplicatively closed subset is also Cohen-Macaulay \cite{k}. As the polynomial ring $A_+ = \mathbb{C}[\sigma_1, \ldots ,\sigma_n]$ is Cohen-Macaulay, so is $A = \mathbb{C}[\sigma_1, \sigma_1^{-1}, \ldots , \sigma_n, \sigma_n^{-1}]$.

\begin{proposition}
For $r \leq n$, the set of points $(a_1, \ldots , a_r) \in A^r$ such that the ideal generated by its coordinates has height $r$, contains an open dense subset of $A^r$. In particular, the set $M$ of points $(a_1, \ldots ,a_n)$ in $A^{n}$ such that the ideal $(a_1, \ldots ,a_n)$ is maximal, contains an open dense subset of $A^{n}$.
\end{proposition}

\noindent Proof: The set of points in $A^r$ corresponding to regular sequences contains an open dense subset of $A^r$. As $A$ is Cohen-Macaulay, the height of an ideal generated by such a sequence equals $r$. The second statement of the proposition now follows because the dimension of $A$ equals $n$.
\hspace*{\fill}$\square$\\

\begin{corr}
Let $r > n$. Then the set of points $(a_1, \ldots ,a_r)$ in $A^r$ such that the ideal $(a_1, \ldots ,a_r)$ equals $A$, contains a Zariski open subset of $A^r$.   
\end{corr}

\noindent Proof: It suffices to prove the statement for $r = n+1$. By the remarks preceding Proposition 3.2, the set of points $a \in M \times A$ such that its coordinates generate a proper ideal is a proper closed subset, hence its complement is open in $M \times A$, and so contains an open dense subset of $A^{n+1}$.  \hspace*{\fill}$\square$\\

\begin{proposition} 
The set of irreducible elements in $A$ contains an open dense subset of the space $A$ when $n > 1$.
\end{proposition}
\noindent Proof: We first prove the statement for the polynomial ring $A_+$. It suffices to show that the complement of the set of irreducible elements in $A_+(d)$ is contained in a proper Zariski closed subset, for all sufficiently large $d$. 

As the degree of a product of two elements in $A_+$ is the sum of the two degrees, every element in $A_+(1)$ is irreducible. 
Now let $d >1$, and let $d_1$, $d_2$, be integers greater than or equal to 1, such that $d_1 + d_2 = d$. Let $\mu : A_+(d_1) \times A_+(d_2) \rightarrow A_+(d)$ be the restriction of the multiplication map in Proposition 3.4. As it is an algebraic map, it follows that $\dim( \overline{\mu (A_+(d_1) \times A_+(d_2)})) \leq \dim (A_+(d_1) \times A_+(d_2))$ (for instance \cite{k}). But $\dim (A_+(d_1) \times A_+(d_2)) = ~{n+d_1 \choose n} + ~{n + d_2 \choose n}$, whereas $\dim (A_+(d_1 + d_2)) = ~{n + d_1 + d_2 \choose n}$. As $n > 1$,  it follows that for $d$ sufficiently large, ${n + d_1 + d_2 \choose n} > ~{n+d_1 \choose n} + ~{n + d_2 \choose n}$ (whereas the reverse inequality is true for $n = 1$); therefore the image of the map $\mu$ is contained in a proper Zariski closed subset of $A_+(d)$.  There are $d-1$ instances of integers $d_1, d_2$ as above which sum to $d$, hence the union of the images of the corresponding multiplication maps is also contained in a proper Zariski closed  subset of $A_+(d)$. This image contains all the elements of degree $d$ that are not irreducible. 

Consider now the ring $A$; being a localisation of a UFD, it is a UFD as well. Hence, the irreducible elements in $A$ correspond to principal prime ideals of $A$. But the prime ideals of $A$ are the prime ideals of $A_+$ that do not intersect the multiplicative closed set ($\sigma_1\sigma_2 ~ \cdots ~\sigma_n$), and the proof now follows as in the last part of Proposition 3.2.
\hspace*{\fill}$\square$\\

\begin{corr} 
In Proposition 3.2, the subset of $P_r$ consisting of elements $(a_1, \ldots ,a_r)$ where each $a_i$ is irreducible also contains an open dense subset of $A^r$ (and similar statements for $U_{r+1}$ and $N_{r+1}$ in Propositions 3.3 and 3.4 respectively). 
\hspace*{\fill}$\square$\\
\end{corr}

We now return to the problem of calculating the degree of autonomy of an $n-D$ system. We recollect that the signal space $\mathbb{C}^{\mathbb{Z}^n}$ is an injective cogenerator \cite{o}, hence $n-D$ behaviors  in $(\mathbb{C}^k)^{\mathbb{Z}^n}$  are in bijective correspondence with submodules of $A^k$. The topology on the set of these behaviors is the Zariski topology on the set of submodules of $A^k$, carried over by the bijection. Genericity statements about submodules then carry over to similar statements about behaviors.\\

Consider an under-determined behavioral system, i.e one defined by a submodule $R \subset A^k$ that can be generated by fewer than $k$ elements (recall from the Introduction that the elements of $R$ are the laws the system obeys). Such a system is non-autonomous, its characteristic ideal is the 0 ideal, and its degree of autonomy equals 0. Thus every behavior in $\mathcal{B}_{\ell,k}, ~ \ell < k$, has degree of autonomy equal to 0 (notation as in Section 2). \\

We now consider over-determined systems. We use the genericity results of Propositions 3.4, 3.5, and Corollary 3.1  to obtain corresponding results for $n-D$ systems.

 Let $\mathcal{M}^*_{\ell,k}$  be the subspace of matrices in $\mathcal{M}_{\ell,k}$ whose characteristic ideals are proper ideals of $A$. By the remark at the beginning of this section, $\frak{i}_{\ell,k}(R) \subseteq \ann(A^k/R) \subseteq \sqrt{\frak{i}_{\ell,k}(R)}$, hence $\frak{i}_{\ell,k}(R) \subsetneq A$ is equivalent to $\ann(A^k/R) \subsetneq A$, which is to say that $R \subsetneq A^k$. Thus $\mathcal{M}^*_{\ell,k}$ is the set of matrices whose $\ell$ rows generate proper submodules of $A^k$.

\begin{thm} 
Let $\ell \geq k$ be such that $\ell - k + 1 \leq n$. Then those matrices which define behaviors whose degree of autonomy equals $\ell - k + 1$ contains an open subset of $\mathcal{M}^*_{\ell,k}$. The other matrices in $\mathcal{M}_{\ell,k}$ define either the zero behavior, or behaviors whose degree of autonomy is strictly less than $\ell - k + 1$.

(ii) Let $\ell - k + 1 > n$. Then an open dense set of matrices in $\mathcal{M}_{\ell,k}$ all define the zero behavior.
\end{thm} 

\noindent Proof: (i) Set $r ={\ell \choose k}$ and $s = \ell - k + 1$.  Recall from Section 2 the commutative diagram (3) and  the continuous map $\frak{m}_{\ell,k}: \mathcal{M}_{\ell,k} \rightarrow A^r$, mapping $R(\sig)$ to its 
$k \times k$ minors (written in some fixed order). By a theorem of Macaulay (for instance, Exercise 10.9 in \cite{k}), the height of the characteristic ideal $\frak{i}_{\ell,k} (R(\sig))$ generated by these minors, is bounded by $s$. We show that this bound is attained by an open  subset of $\mathcal{M}^*_{\ell,k}$.

Consider the following matrix in $\mathcal{M}_{\ell,k}$

\[ R(\sig)  = \begin{pmatrix} a_1 & 0 & \cdots & 0 \\
0 &  & \\  
\vdots & &  I_{k-1} \\
0 &  &  \\ 
a_2  & 0 & \cdots & 0 \\ 
\vdots & & \\
a_s & 0 & \cdots & 0
\end{pmatrix} 
\]
where the $(k-1) \times (k-1)$ identity matrix $I_{k-1}$ occupies rows and columns 2 to $k$, and the $a_i$ are nonzero elements of $A$. Its nonzero maximal minors are 
$\{a_1, \ldots , a_s \}$. 
If the sequence defined by these minors were regular, then the depth of the ideal $\frak{i}_{\ell,k}(R(\sig))$  generated by  them, and hence its height, would be at least $s$.  As the number of generators is exactly $s$, its height would be equal to $s$ (by Krull's Height Theorem).

The  entries $\{a_i\}$ of the above $R(\sig)$ could be arbitrary elements from $A$, hence it follows that the composition  $\mathcal{M}_{\ell,k} \stackrel {\frak{m}_{\ell,k}} {\longrightarrow} A^r \stackrel{\pi}{\longrightarrow} A^s$ is surjective,  where the second map is the projection to  $A^s$ (the indices of $s$ determined by the nonzero minors described above).
As $s \leq n$, the set of points in $A^s$ whose coordinates generate ideals of $A$ of
height $s$ contains an open dense subset (Proposition 3.5), hence it follows that there is a nonempty open dense
subset of matrices in $\mathcal{M}_{\ell,k}$ with the property that the ideal generated by the $s$ minors described above is of height $s$. Intersecting this open set with $\mathcal{M}^*_{\ell,k}$ gives the open subset of matrices of the theorem whose characteritic ideals are proper,  have height at least $s$, and hence height equal to $s$ by the theorem of Macaulay quoted above. 
These matrices define behaviors whose degree of autonomy equals $s$ (by Wood, Rogers and Owens
\cite{wo} and Napp and Rocha \cite{d}), and this is the maximum possible degree of autonomy of a nonzero behavior defined by matrices in $\mathcal{M}_{\ell,k}$.

(ii) If $s = \ell - k + 1 > n$, then by Corollary 3.1 there is an open dense subset of points in $A^s$ whose coordinates generate the unit ideal in $A$. Let $(a_1, \ldots , a_s)$ be such a point, then there is a matrix just as in (i) above, whose maximal minors are these $a_i$.  Hence the set $U$ of matrices in
 $\mathcal{M}_{\ell,k}$ such that for every $R(\sig) \in U$, the ideal generated by the corresponding $s$ minors
equals $A$, is nonempty open, and so open dense. Thus $\frak{i}_{\ell,k}(R(\sig)) = A$ and the rows of $R(\sig)$ generate $A$, for every $R(\sig) \in U$. All these matrices then define the zero behavior. \hspace*{\fill}$\square$\\

\noindent Remark: Suppose $n$, $\ell$, and $k$, $\ell \geq k$, are such that $\mathcal{M}^*_{\ell,k}$ is open in $\mathcal{M}_{\ell,k}$; then the open subset of $\mathcal{M}^*_{\ell,k}$ described in the above theorem would be  open and dense in $\mathcal{M}_{\ell,k}$. For example, if $k = 1$, then $\ell \leq n$ is such an instance by Proposition 3.2. Another instance is $k = \ell$ when there is only one maximal minor, then the set $\mathcal{M}^*_{\ell,\ell}$ of square matrices whose determinants are non-units is open dense in $\mathcal{M}_{\ell,\ell}$ (independent of $n \geq 1$, by the remarks preceeding Lemma 2.2).

If $\ell < k$, then every $R(\sig) \in \mathcal{M}_{\ell,k}$ satisfies $\frak{i}_{\ell,k}(R(\sig)) \subsetneq A$; indeed the characteristic ideal is equal to 0. \\

We can now prove the main result of the section.  Let $\mathcal{S}^*_{\ell,k} =  \mathcal{S}_{\ell,k} \setminus \{A^k\}$ be the subspace of proper submodules of $A^k$ that can be generated by $\ell$ elements, and let $\mathcal{I}^*_t = \mathcal{I}_t \setminus \{A\}$ be the subspace of proper ideals that can be generated by $t$ elements (notation as in Section 2). By the observation preceeding Theorem 3.1, the maps $\Pi_{\ell,k}$ and  $\frak{i}_{\ell,k}$ of diagram (3) restrict to   $\Pi_{\ell,k}: \mathcal{M}^*_{\ell,k} \rightarrow \mathcal{S}^*_{\ell,k}$ and $\frak{i}_{\ell,k}: \mathcal{S}^*_{\ell,k} \rightarrow \mathcal{I}^*_{\ell \choose k}$. Let $\mathcal{B}^*_{\ell,k} = \mathcal{B}_{\ell,k} \setminus \{ 0\}$ be the subspace of nonzero behaviors in bijective correspondence with $\mathcal{S}^*_{\ell,k}$. 
    
\begin{thm}
(i) Let $\ell \geq k$ be such that  $\ell - k + 1 \leq n$. Then the set of behaviors in 
$\mathcal{B}^*_{\ell,k}$ whose degree of autonomy equals $\ell - k + 1$, contains an open 
subset of $\mathcal{B}^*_{\ell,k}$. The other behaviors have degree of autonomy strictly less than than $\ell - k + 1$.

(ii) If  $\ell - k + 1 > n$, then the zero behavior is open dense in $\mathcal{B}_{\ell,k}$.
Thus, this open dense behavior has degree of autonomy equal to $\infty$.
\end{thm}

\noindent Proof: (i) Set $r = {\ell \choose k}$ and $s =~\ell - k + 1$. In the above notation,
$\mathcal{I}^*_s \subset \mathcal{I}_s$ is the subspace of proper ideals that can be generated by $s$ elements.  Let
$X \subset \mathcal{I}^*_s$ be the subset of those ideals whose heights equal $s$.
We show $X$ contains an open subset of $\mathcal{I}^*_s$, and as $\mathcal{I}^*_s$ is equipped
with the quotient topology given by the surjection $\Pi_{s,1}: P_s \rightarrow \mathcal{I}^*_s$, we
need to show that $U := \Pi_{s,1}^{-1}(X)$ contains an open subset of $P_s$ ($P_s$ as in Proposition 3.2).

Let $I \in X$, and let $(a_1, \ldots ,a_s) \in \Pi_{s,1}^{-1}(I)$. If $a_1, \ldots ,a_s$ is a
regular sequence, then there is nothing to be done by Proposition 3.4. Suppose that it is not.
As $A$ is Cohen-Macaulay, $\dep(I) = s$, hence there is a regular sequence $a'_1, \ldots ,a'_s$
in $I$.
Each $a'_i$ is an $A$-linear combination of $a_1, \ldots ,a_s$, and it follows that there is an
$A$-linear map $L: A^s \rightarrow A^s$ such that $L(a_1, \ldots ,a_s) = (a'_1, \ldots ,a'_s)$.
By Proposition 3.4, there is a neighbourhood $W$ of $(a'_1, \ldots ,a'_s)$ containing an open set
of points, each of whose coordinates is a regular sequence. The map $L$ is continuous, hence
$L^{-1}(W)$ contains an open neighbourhood of $(a_1, \ldots ,a_s)$, and let $(b_1, \ldots ,b_s)$ be
a point in it. Then $L(b_1, \ldots ,b_s) = (b'_1, \ldots ,b'_s)$ is in $W$, and so the sequence
$b'_1, \ldots ,b'_s$ is regular. Each $b'_i$ belongs to the ideal $J = (b_1, \ldots , b_s)$,
hence $\dep(J) = s$, and so also $\height(J) = s$. Thus  $L^{-1}(W) \subset U$, and $U$ is open.

It now follows that the subset of ideals in $\mathcal{I}^*_r$, which project to the collection $X$ under the map $\pi: A^r \rightarrow A^s$ (in the proof of the preceding theorem), is open.  Denote this subset by $X_r$.
Then $X' = \frak{i}_{\ell,k}^{-1}(X_r)$ is open in $\mathcal{S}^*_{\ell,k}$,
and the subset of behaviors in $\mathcal{B}^*_{\ell,k}$ corresponding to submodules
in $X'$ is also open. All their degrees of autonomy  equal $\ell - k + 1$, and this is the maximum possible by Macaulay's Theorem.

(ii) The zero behavior in $\mathcal{B}_{\ell,k}$ corresponds to $A^k$ in $\mathcal{S}_{\ell,k}$.
By (ii) of Theorem 3.1, there is an open dense set of matrices in $\mathcal{M}_{\ell,k}$ whose
rows generate $A^k$. As the topology of $\mathcal{S}_{\ell,k}$ is given by the surjection
$\Pi_{\ell,k}: \mathcal{M}_{\ell,k} \rightarrow \mathcal{S}_{\ell,k}$ (diagram (3) of Section 2),
it follows that the zero behavior is open and dense in $\mathcal{B}_{\ell,k}$.
\hspace*{\fill}$\square$

\section{Strength of controllers}

We use the results of the above section to study notions of {\em strength} and {\em efficiency} of a controller. \\

The notion of autonomy that we have studied is intimately related to the following questions: How are the trajectories of an $n-D$ system determined? Is there a sub-lattice of $\mathbb{Z}^n$ such that these trajectories assume arbitrary values on it? Do the values of a trajectory on the points of a proper sub-lattice determine it on all of $\mathbb{Z}^n$? 

For example, consider the $2-D$ system defined by the kernel of the map
\[ \left(\begin{array}{c} \sigma_1 - 1 \\ \sigma_2 - 1 \end{array} \right): (\mathbb{C})^{\mathbb{Z}^2} \rightarrow (\mathbb{C}^2)^{\mathbb{Z}^2}
\]
The trajectories of this system are the constants, and a trajectory is determined by its value at $0 \in \mathbb{Z}^2$. This is captured by the fact that the degree of autonomy of this system equals 2, the codimension of the 0 sub-lattice. On the other hand, if we consider the system whose trajectories are the kernel of
\[ \sigma_1-1: (\mathbb{C})^{\mathbb{Z}^2} \rightarrow (\mathbb{C})^{\mathbb{Z}^2} \] then values arbitrarily assigned on the points $(0, n) \in \mathbb{Z}^2$ results in a unique trajectory, and the degree of autonomy of this system is 1. The degree of autonomy of a system measures this possibility of assigning initial values on a sub-lattice of $\mathbb{Z}^n$.

When a behavior must be restricted to a sub-behavior satisfying certain properties (of stability, of growth at infinity, and so on) by the attachment of a controller, we can generically include the condition that trajectories of the controlled system be determined by initial values  prescribed on as small a sub-lattice as possible. This is the content of the results of this paper. 

We can isolate this particular property of reducing the dimension of the sub-lattice on which initial values be prescribed, and consider a controller {\em efficient} if it can reduce the dimension maximally amongst all controllers with the same number of laws (equations). This leads us to the following considerations.\\

Define $\delta: \mathcal{B}(k) \rightarrow \{0, 1, \ldots ,n\} \cup \{ \infty \}$, mapping a behavior $B$ in $(\mathbb{C}^k)^{\mathbb{Z}^n}$ to $\delta (B)$, its degree of autonomy. 

\vspace{.1cm}

Suppose $B = \mathcal{B}(R)$ is an $n-D$ behavior defined by a submodule $R$ of $A^k$. Suppose we attach a controller $C = \mathcal{B}(R')$ to it, to obtain the controlled system $B \cap C = \mathcal{B}(R + R')$. 

\vspace{.2cm}

Define $\varsigma: \mathcal{B}(k) \times \mathcal{B}(k) \rightarrow \{0, 1, 2, \ldots ,n\} \cup \{ \infty \}$, mapping $(C,B)$ to $\varsigma(C,B) = \delta (B \cap C) - \delta (B)$ for $B \neq 0$, and $\varsigma(C,0) = \infty$ for every $C$. 

\vspace{.1cm}

$\varsigma(C,B)$ is the {\em strength} of $C$ with respect to $B$, it measures the increase in the degree of autonomy of the behavior $B$ upon attaching the controller $C$ to it. 

\vspace{.1cm}

The controller $C$ is said to be {\em maximally efficient} with respect to $B$, if for all $C'$ with 
$\delta (C') = \delta (C)$, $\varsigma(C',B) \leq \varsigma(C,B)$. \\

\begin{thm} Let $\ell < k$, and let $B$ be in $\mathcal{B}_{\ell,k} \setminus \mathcal{B}_{\ell -1,k}$
(thus $B$ is a non-autonomous system).

(i) Let $\ell '$ be such that  $1 \leq \ell + \ell ' - k + 1 \leq n$. Then there is an open subset
$U$ in the space of nonzero behaviors $\mathcal{B}^*_{\ell',k}$ such that for every $C \in U$, the strength of $C$ with
respect to $B$ is $\varsigma(C,B) =~ \ell + \ell ' - k + 1$. Every such $C$ is thus maximally
efficient with respect to $B$. If $\ell' \leq k$, then $U$ is open dense in $\mathcal{B}_{\ell',k}$.

(ii) If $\ell + \ell ' - k + 1 > n$, then there is an open dense $U \subset \mathcal{B}_{\ell',k}$
such that for every $C$ in it, the controlled system $B \cap C$ is the zero behavior, and hence $\varsigma(C,B) = \infty$.
\end{thm}

\noindent Proof: (i) Set $s = ~ \ell + \ell ' - k + 1$. Let
$B \in \mathcal{B}_{\ell,k} \setminus \mathcal{B}_{\ell - 1,k}$ be defined by the submodule
$R \subset A^k$, where $R$ can be generated by $\ell$ elements but not  by $\ell - 1$ elements.
Let $M$ be an $\ell \times k$ matrix whose rows generate $R$.
Define $\frak{m}_R: \mathcal{M}^*_{\ell', k} \rightarrow 
A^s$ by mapping a matrix $M'$ with $\ell'$ rows, to the $k \times k$ minors
described in Theorem 3.1 of the matrix obtained by adjoining the rows of $M'$ to the rows of $M$. This is a continuous map,
hence the inverse image of the open set of points in $A^s$  whose coordinates define ideals of
height $s$ (Proposition 3.5), is open in $\mathcal{M}^*_{\ell',k}$. By restricting to a smaller open set if necessary, we can ensure that the matrix obtained by adjoining the rows of $M'$ to $M$ lies in $\mathcal{M}^*_{\ell + \ell'}$. This map descends to $\mathcal{S}^*_{\ell',k} \rightarrow \mathcal{I}^*_s$ to define an open set $U$ in $\mathcal{B}^*_{\ell',k}$ such that for every $C$ in $U$, $\delta (B \cap C) = s$. As $\delta (B) = 0$, $\varsigma (C,B) = s$.

A behavior in $\mathcal{B}^*_{\ell + \ell',k}$ has degree of autonomy at most $s$, hence every
controller $C$ in $U$ is maximally efficient with respect to $B$.

The last claim follows from the remark following Theorem 3.1.

(ii) The proof  follows similarly by an application of Theorem 3.1 (ii), and Corollary 3.1 instead of Proposition 3.5.
\hspace*{\fill}$\square$\\

\begin{thm}
(i) Let $\ell \geq k$ and $\ell'$ be such that $1 \leq \ell + \ell' - k + 1 \leq n$. Then there are open
subsets $U$ of $\mathcal{B}^*_{\ell,k}$ and $U'$ of $\mathcal{B}^*_{\ell',k}$ such that
for every $B \in U$ and $C \in U'$, $\varsigma(C,B) = ~\ell'$. Thus every controller in $U'$
is maximally efficient with respect to every $B$ in $U$.

(ii) If $\ell$ and $\ell'$ are such that $\ell + \ell ' - k + 1  > n$, then $B \cap C = 0$
and $\varsigma(C,B) = \infty$ for all $B$ and $C$ belonging to open neighbourhoods in
$\mathcal{B}_{\ell,k}$ and $\mathcal{B}_{\ell',k}$.
\end{thm}

\noindent Proof:  (i) Consider the map $\alpha: \mathcal{M}_{\ell,k} \times \mathcal{M}_{\ell',k} 
\rightarrow \mathcal{M}_{\ell + \ell',k}$ given by adjoining the $\ell'$ rows of one matrix to
the $\ell$ rows of the other. It is continuous, and surjective.
Let $W$ be the open subset of $\mathcal{M}^*_{\ell + \ell',k}$ consisting of those matrices
that define behaviors whose degree of autonomy equals $\ell + \ell' - k + 1$ (Theorem 3.1).
Let $V$ and $V'$ be the projections of  $\alpha ^{-1} (W)$ to the two factors. Their intersections with open subsets of matrices in
$\mathcal{M}^*_{\ell,k}$ and  $\mathcal{M}^*_{\ell',k}$ defining behaviors of maximum degree of
autonomy, give the open neighbourhoods $U$ and $U'$ of the theorem. The proof of (ii) is the same as the proof of 4.1 (ii) above.
\hspace*{\fill}$\square$\\

\noindent Remark: The notion of efficiency we have proposed in this paper is a reflection of the idea that an efficient controller should have little in common with the the system it regulates. A related notion introduced by Willems is that of a `regular interconnection' \cite{w1}. The relationship between these two notions is subtle, and will be pursued elsewhere. \\

\noindent Remark: Pal and Pillai \cite{pp} have established results on degree of autonomy, similar to those of \cite{d,wo}, for {\em scalar} distributed systems defined by PDE, i.e. systems defined by ideals of $\mathbb{C}[\partial_1, \ldots, \partial_n]$ in the space of smooth functions ($k = 1$ in the notation of the paper). All the results of this paper on genericity hold in their situation as well. Such results are not available for general systems of PDE. If they were, then the genericity results of this paper would also follow in this  case. 

\vspace{.5cm}
 
\section{Acknowledgement} We are extremely grateful to Bharat Adsul, Krishna Hanumanthu, Debasattam
Pal, Ananth Shankar,  Arul Shankar and Jugal Verma for many useful conversations. We are grateful to the referees
for their comments and for saving us from error. The first author is also grateful to the Department
of Electrical and Computer Engineering, University of Porto, for its hospitality.

\end{document}